\documentclass{esannV2}
\usepackage[dvips]{graphicx}
\usepackage[latin1]{inputenc}
\usepackage{amssymb,amsmath,array}
\usepackage[textwidth=2.5cm]{todonotes}
\usepackage[ruled,norelsize]{algorithm2e}
\usepackage{comment}
\newcommand{\pierre}[1]{\todo[color=purple]{{\bf Pierre:} #1}}
\newtheorem{lemma}{Lemma}

\begin{document}
%style file for ESANN manuscripts
\title{Beyond Pham's algorithm for joint diagonalization}

%***********************************************************************
% AUTHORS INFORMATION AREA
%***********************************************************************
\author{Pierre Ablin$^1$, Jean-Fran\c cois Cardoso$^2$ and Alexandre Gramfort$^1$
%
% Optional short acknowledgment: remove next line if non-needed
\thanks{This work was supported by the Center for Data Science, funded by the IDEX Paris-Saclay, ANR-11-IDEX-0003-02, and the European Research Council (ERC SLAB-YStG-676943).}
%
% DO NOT MODIFY THE FOLLOWING '\vspace' ARGUMENT
\vspace{.3cm}\\
%
% Addresses and institutions (remove "1- " in case of a single institution)
1- INRIA - Parietal team \\
1 Rue Honor\'e d'Estienne d'Orves, 91120 Palaiseau - France
%
% Remove the next three lines in case of a single institution
\vspace{.1cm}\\
2- CNRS - Institut d'Astrophysique de Paris \\
98bis boulevard Arago, 75014, Paris - France\\
}
%***********************************************************************
% END OF AUTHORS INFORMATION AREA
%***********************************************************************

\maketitle

\begin{abstract}
The approximate joint diagonalization of a set of matrices consists in finding a basis in which these matrices are as diagonal as possible.
This problem naturally appears in several statistical learning tasks such as blind signal separation.
We consider the diagonalization criterion studied in a seminal paper by Pham (2001), and propose a new quasi-Newton method for its optimization.
Through numerical experiments on simulated and real datasets, we show that the proposed method outperforms Pham's algorithm.
An open source Python package is released.
\end{abstract}

\section{Introduction}
\label{sec:intro}
The task of joint diagonalization arises in several formulations of the blind source separation problem.
In~\cite{cardoso1993blind}, independent component analysis is performed by joint-diagonalization of a set of cumulant matrices.
In~\cite{belouchrani1997blind}, the separation of non-stationnary signals is carried by joint-diagonalization of a set of autocorrelation matrices.
Finally, in~\cite{pham2001blind}, joint-diagonalization of a set of covariance matrices separates Gaussian sources that have non-stationnary power.

Consider a set of $n$ symmetric square matrices $\mathbf{\mathcal{C}} \triangleq (C^1,\cdots, C^n)$ of size $p \times p$.
Its approximate joint diagonalization consists in finding a matrix $B \in \mathbb{R}^{p \times p}$ such that the matrix set $B \mathbf{\mathcal{C}}B^{\top} \triangleq (B C^1B^{\top}, \ldots, B C^n B^{\top})$ contains matrices that are \emph{as diagonal as possible}, as measured by some joint-diagonality criterion.
This paper considers the joint diagonalization of positive matrices,
defined as the minimization of the (non-convex) criterion
\vspace{-0.2cm}
\begin{equation}
\label{eq:loss}
    \mathcal{L}(B) =\frac{1}{2n} \sum_{i=1}^n \Big [\log\det \text{diag}(BC^iB^{\top}) - \log\det(BC^iB^{\top}) \Big ]  \enspace, \vspace{-0.3cm}
\end{equation}
This criterion was introduced  by Pham in~\cite{pham2001blind}
who derived it as the negative log-likelihood of a source separation model for
Gaussian stationary sources.

\begin{comment}
Other criteria have been proposed in the literature. \pierre{Maybe cut this part if the paper is too long}
%
First, one can use the Frobenius criterion~\cite{cardoso1993blind}, which writes: $\mathcal{L}(B) = \sum_{i=1}^n ||\text{OFF}(BC^iB^{\top}||_F^2$, where the $\text{OFF}$ of a matrix are the off-diagonal coefficients, and $||\cdot||_F$ is the Frobenius norm.
%
This criterion is natural, but suffers from two drawbacks.
%
First, it does not stem from a statistical framework and is rather ad-hoc.
%
Second, it is trivially minimized by taking $B=0$.
%
To overcome this problem, it is common to constrain $B$ to orthogonality~\cite{cardoso1993blind}.
%
The set can also be directly diagonalized in the least squares sense~\cite{yeredor2002non} by minimizing $\sum_{i=1}^n||C^i - B^{\top} \Lambda_i B||_F^2$ with respect to $B$ and the diagonal matrices $\Lambda_1, \cdots, \Lambda_n$.
%
\end{comment}

%
Pham~\cite{pham2001joint} proposes in its seminal work a block coordinate descent approach for the minimization of $\mathcal{L}$.
Each iteration of this method guarantees a decrease of $\mathcal{L}$.
Further, when there exists a matrix $B_*$ such that $B_*\mathcal{C}B_*^{\top}$ contains only diagonal matrices (that is, if the set is exactly jointly diagonalizable), then in the vicinity of $B_*$, the algorithm converges quadratically.
Since then, only a few other papers have focused on minimizing this criterion.
In~\cite{joho2008newton}, the Newton method is studied.
However, this method is not practical as it requires solving a $p^2\times p^2$ linear system at each iteration.
In~\cite{todros2007fast}, it is proposed to minimize an approximation of the cost function.
The authors do so by alternating minimization over columns.
This method minimizes an approximation of $\mathcal{L}$, and scales in $O(p^4)$.

In this paper, we propose a quasi-Newton method for the minimization of $\mathcal{L}$.
We use sparse approximations of the Hessian which are cheap to compute and match the true Hessian when the set is jointly diagonalized, granting quadratic convergence.
The algorithm as a cost per iteration of $O(n \times p^2)$, which is the natural scaling of the problem since it is the size of the dataset.
Through experiments, we show that the proposed method outperforms Pham's algorithm, on both synthetic and real data.

\textbf{Notation}: The identity matrix of size $p$ is denoted $I_p$.
The Frobenius scalar product between two $p\times p$ matrices is noted as $\langle M |M' \rangle \triangleq \sum_{a, b} M_{ab}M'_{ab}$.
The corresponding norm is $||M||_F = \sqrt{\langle M |M \rangle }$.
Given a $p\times p \times p \times p$ tensor $\mathcal{H}$, the weighted scalar product is $\langle M |\mathcal{H} |M' \rangle \triangleq \sum_{a, b, c, d} \mathcal{H}_{abcd}M_{ab}M'_{cd}$.
The Kronecker symbol $\delta_{ab}$ is $1$ is $a=b$, $0$ otherwise.
\section{Study of the cost function}
\label{sec:loss}

The cost function $\mathcal{L}(B)$ is defined on the group of invertible matrices,
with $\log\det(BC^iB^{\top})$ acting as a barrier.
%
%%Hence, the problem is non-convex.
%
Its minimization is performed by iterative algorithms.
To exploit the multiplicative structure of the group of invertible matrices,
we perform \emph{relative} updates on $B$~\cite{cardoso1996equivariant}.
The neighborhood of an iterate $B^{(t)}$ is parameterized by a small $p\times p$ matrix $\mathcal{E}$ as  $B^{(t+1)} = (I_p + \mathcal{E})B^{(t)}$.
The second-order (in $\mathcal{E}$) Taylor expansion of the loss function 
\vspace{-0.2cm}
\begin{equation}
    \label{eq:taylor}
    \mathcal{L}((I_p + \mathcal{E})B) = \mathcal{L}(B) + \langle G |\mathcal{E}\rangle + \frac12 \langle \mathcal{E} | \mathcal{H} | \mathcal{E} \rangle + o(||\mathcal{E}||^2) \enspace .\vspace{-0.3cm}
\end{equation}
where the $p \times p$ matrix $G(B)$ is the relative gradient
and the $p \times p \times p \times p$ tensor $\mathcal{H}(B)$
is the relative Hessian.
Simple algebra yields:
\vspace{-0.2cm}
\begin{equation}
    \label{eq:hessian}
    G_{ab} = \frac{1}{n} \sum_{i=1}^n \frac{D^i_{ab}}{D^i_{aa}} - \delta_{ab}, \enspace \text{and}  \enspace
    \mathcal{H}_{abcd} = \delta_{ac} \frac{1}{n} \sum_{i=1}^n \left[\frac{D^i_{bd}}{D^i_{aa}}-2\frac{D^i_{ab}D^i_{ad}}{{D^i_{aa}}^2}\right]+ \delta_{ad}\delta_{bc} \enspace
\end{equation}
where we define $D^i = BC^iB^{\top}$ for $i=1, \dots, n$.
Eq.~\eqref{eq:hessian} shows that $H_{aaaa}=0$,
consistent with the scale-indeterminacy of the criterion: %
$\mathcal{L}(\Lambda B)=\mathcal{L}(B)$ for any diagonal matrix $\Lambda$. 
The criterion is flat in the direction of the scale matrices.

\textbf{Complexity}: The cost of computing a gradient is $O(p^2 \times n)$.
It is the natural complexity of an iterative algorithm as it is the same cost as computing the set $B\mathcal{C}B^{\top}$.
Computing the Hessian is $O(p^3 \times n)$, and computing $H^{-1}G$ is $O(p^6)$.
This is prohibitively costly when $p$ is large compared to a gradient evaluation.
As a consequence, Newton's method is extremely costly to setup.

\section{Relative Quasi-Newton method}
\label{sec:quasinewton}

In this section, we first introduce an approximation of the Hessian and then derive a quasi-Newton algorithm to minimize~\eqref{eq:loss}.
\subsection{Hessian approximation}

The accelerated algorithm presented in this article is based on the massive
sparsification of the Hessian tensor when matrices $D^i$ are all diagonal.
Indeed, it that case, it becomes:
\vspace{-0.2cm}
\begin{equation}
    \label{eq:approx}
     \tilde{\mathcal{H}}_{abcd} = 
     \delta_{ac}\delta_{bd} \Gamma_{ab} + \delta_{ad}\delta_{bc} - 2 \delta_{abcd}
     \qquad
     \Gamma_{ab} =\frac{1}{n} \sum_{i=1}^n \frac{D^i_{bb}}{D^i_{aa}}
\end{equation}
%
%Define $\Gamma_{ab} =\frac{1}{n} \sum_{i=1}^n \frac{D^i_{bb}}{D^i_{aa}}$ such that %$\tilde{\mathcal{H}}_{abcd} = \delta_{ac}\delta_{bd} \Gamma_{ab} + %\delta_{ad}\delta_{bc} - 2 \delta_{abcd}$.
%
This Hessian approximation has three key properties.
First, it is cheap to compute, at cost $O(p^2 \times n)$, just like a gradient.
Then, it is sparse and structured. It only has $\simeq 2p^2$ non-zero coefficients, and can be seen as a block-diagonal operator, with blocks of size $2$.
Indeed, for a $p\times p$ matrix $M$, we have: 
\begin{equation}
    \label{eq:block}
    \begin{pmatrix}
        [\tilde{\mathcal{H}}M]_{ab}      \\
        [\tilde{\mathcal{H}}M]_{ba}       
    \end{pmatrix}
    = 
    \begin{pmatrix}
        \Gamma_{ab} & 1 - 2\delta_{ab} \\
        1 - 2\delta_{ab}  & \Gamma_{ba}       
    \end{pmatrix} 
    \begin{pmatrix}
        M_{ab}      \\
        M_{ba}       
    \end{pmatrix}
    \enspace .
\end{equation}
The following lemma establishes the positivity of the approximate Hessian:

\begin{lemma}
\label{lemma:pos}
The approximation $\tilde{\mathcal{H}}$ is positive with $p$ zero eigenvalues.
If the matrices $C^i$ are independently sampled from continuous densities,  with probability one, the other $p^2 - p$ eigenvalues are strictly positive.
\end{lemma}

\textbf{Proof:} Using eq.~\eqref{eq:approx}, one has $\tilde{\mathcal{H}}E_{ii} = 0$, where $E_{ii}$ is the matrix with $1$ for its $(i, i)$ coefficient, and $0$ elsewhere.
Thus $\tilde{\mathcal{H}}$ has $p$ zero eigenvalues, and the associated eigenvectors are the $E_{ii}$ for $i =1\cdots p$.
The $p^2 - p$ other eigenvalues are the eigenvalues of the $2\times 2$ blocks introduced in eq.~\eqref{eq:block}.
The diagonal coefficients of the blocks are the $\Gamma_{ab} > 0$.
The determinant of a block is given by $\Gamma_{ab}\Gamma_{ba} - 1$.
Cauchy-Schwartz inequality yields $\Gamma_{ab}\Gamma_{ba} \geq 1$, with equality if and only if $D^i_{aa} \propto D^i_{bb}$.
This happens with probability $0$, concluding the proof.

Finally, the approximation is straightforwardly inverted by inverting each $2 \times 2$ blocks.
The Moore-Penrose pseudoinverse of $\tilde{\mathcal{H}}$, $\tilde{\mathcal{H}}^+$, satisfies:
\vspace{-0.2cm}
\begin{equation}
    \label{eq:inverse}
    [\tilde{\mathcal{H}}^+G]_{ab} = \frac{\Gamma_{ba}G_{ab} - G_{ba}}{\Gamma_{ab}\Gamma_{ba} - 1}, \enspace \forall a \neq b \enspace, \vspace{-0.3cm}
\end{equation}
and $[\tilde{\mathcal{H}}^+G]_{aa} = 0$.
The cost of inversion is thus $O(p^2)$.

\subsection{Algorithm}
\label{sec:algo}

\begin{algorithm}[tb]
\SetKwInOut{Input}{Input}
\SetKwInOut{Output}{Output}
 \Input{Set of matrices $\mathcal{C}$, number of iterations $T$.}
 Initialize $B=I_p$. \\
 \For{$t=1,\cdots, T$}{
    Compute the gradient $G$ using eq.\eqref{eq:hessian} \\
    Compute the Hessian approximation $\tilde{\mathcal{H}}$ using eq.\ref{eq:approx}\\
    Compute the search direction $-\tilde{\mathcal{H}}^{+}G$ using eq.~\eqref{eq:inverse} \\
    Do a backtracking line search in that direction to find a step size $\alpha$ decreasing $\mathcal{L}$\\
    Set $B \leftarrow (I_p - \alpha\tilde{\mathcal{H}}^{+}G) B $ \\
  }
 \Output{ $B$}
 \caption{Quasi-Newton method for joint-diagonalization}
 \label{algo:qn}
\end{algorithm}
The proposed quasi-Newton method uses $-\tilde{\mathcal{H}}^{+}G$ as search direction.
Following from the positivity of $\tilde{\mathcal{H}}$, this is a descent direction.
Unfortunately, there is no guarantee that the iteration $B\leftarrow (I_p -\tilde{\mathcal{H}}^{+}G)B$ decreases the cost function.
Therefore, we resort to a line-search to find a step $\alpha >0$ ensuring $\mathcal{L}((I_p -\alpha\tilde{\mathcal{H}}^{+}G)B)< \mathcal{L}(B)$ (condition $(*)$).
This is done using backtracking, starting from $\alpha = 1$ and iterating $\alpha \leftarrow \frac{\alpha}{2}$ until the condition $(*)$ is met.
The full algorithm is summarized in Algorithm.~\ref{algo:qn}.

\textbf{Quadratic convergence}: Like Pham's algorithm, the proposed algorithm enjoys quadratic convergence when the matrix  set is jointly diagonalizable.
Indeed, assume that there exists a matrix $B_*$ such that $B_*\mathcal{C}B_*^{\top}$ contains only diagonal matrices.
Then, by construction, $\mathcal{H}(B_*) = \tilde{\mathcal{H}}(B_*)$.
It follows that the method converges quadratically in the vicinity of $B_*$.

\section{Experiments}
\label{sec:expe}

\subsection{Setting}
For the experiments, three data sets are used -- coming either from synthetic or real data -- and we compare our approach to Pham's algorithm.
The code to reproduce the experiments is available online\footnote{https://github.com/pierreablin/qndiag}.
We set $n=100$ and $p=40$.

\textbf{Initialization}: For a dataset $\mathcal{C}$, the algorithms start from a $\emph{whitener}$ (whitening matrix):  writing $PDP^{\top}=\frac{1}{n}\sum_{i=1}^n C^i$ with $P$ orthogonal and $D$ diagonal, the initial matrix is taken as $B^{(0)} = D^{-\frac12}P^{\top}$.

\textbf{Metrics}: To compare the speed of convergence of the algorithms, we monitor the diagonalization error $\mathcal{L}(B)$, and the gradient norm $||G(B)||$.
The first quantity goes to $0$ if the dataset is perfectly diagonalizable, while the second should always converge to $0$ since the algorithm should reach a local minimum.

\textbf{Synthetic data}: We proceed as in~\cite{ziehe2004fast} for generating synthetic datasets.
We generate $n$ diagonal matrices of size $p \times p$, $(D^1,\cdots D^n)$ for which each diagonal coefficient is drawn from an uniform distribution in $[0, 1]$.
Then, we generate a random `mixing' matrix $A\in \mathbb{R}^{p \times p}$ with independent normally distributed entries.
Finally, in order to add noise, we generate $n$ matrices $R^1, \cdots, R^n \in \mathbb{R}^{p \times p}$ of normally distributed entries.
The dataset is then $\mathcal{C} = (C^1,\cdots, C^n)$ with:
\vspace{-0.2cm}
\begin{equation}
    C^i = A D ^i A ^{\top} + \sigma^2 (R^i) (R^i)^{\top} \enspace, \vspace{-0.3cm}
\end{equation}
where $\sigma$ controls the noise level.
In practice we take $\sigma=0$ (experiment \textbf{(a)}, perfecty jointly-diagonalizable set) or $\sigma = 0.1$ (experiment \textbf{(b)}).

\begin{figure}
    \centering
    \includegraphics[width=\columnwidth]{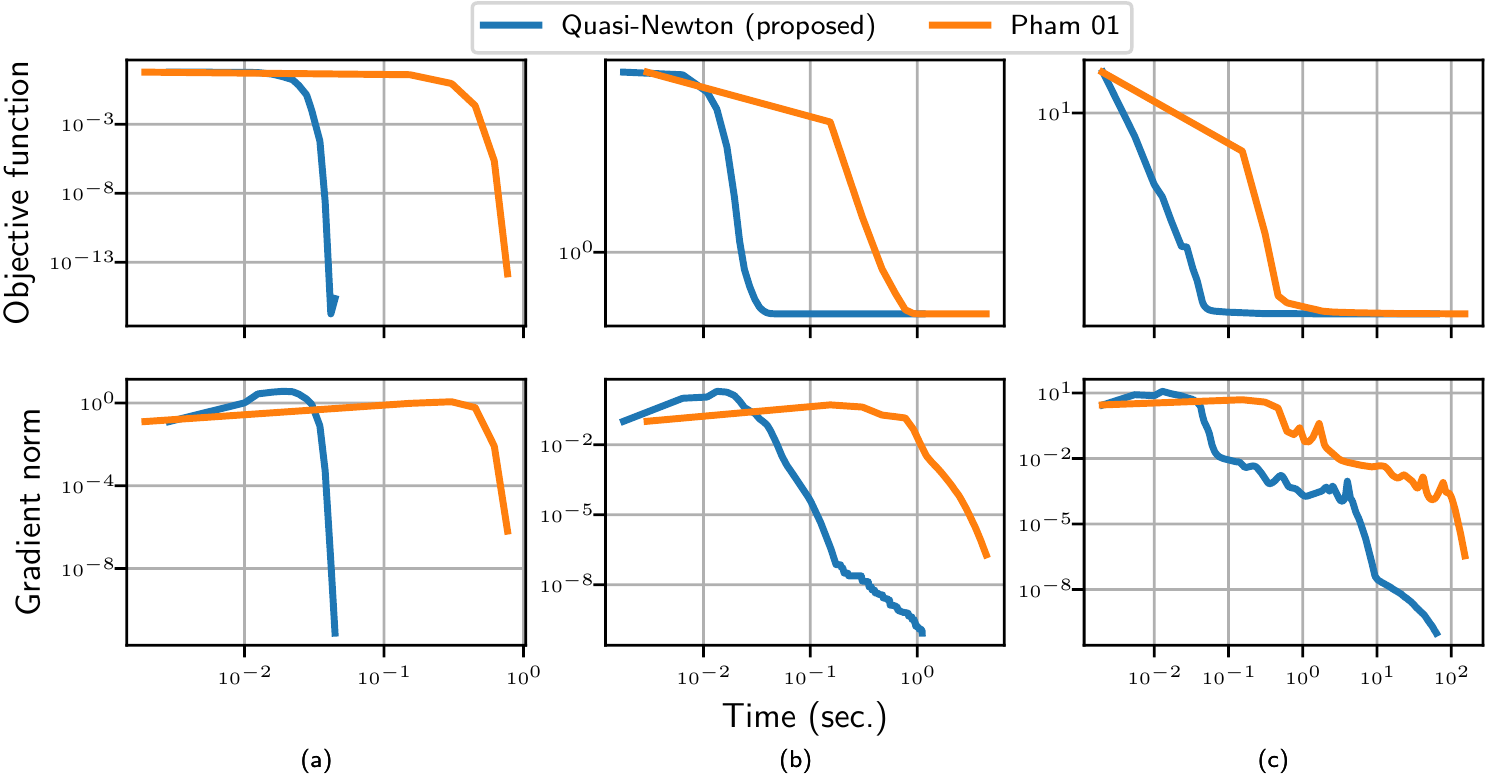}
    \caption{Comparison of the algorithms on three datasets. \textbf{(a)}: Jointly diagonalizable synthetic dataset. \textbf{(b)}: Same dataset with added noise, rendering perfect joint diagonalization impossible.\textbf{(c)}: Real data, covariance matrices from MEG signals. Note the log-log scale.}
    \label{fig:results}
\end{figure}

\textbf{Magnetoencephalography (MEG) data}: We use the MNE sample dataset~\cite{gramfort2014mne}.  From $n$ matrices containing $p$ signals of $T$ samples, $X_1, \cdots, X_n \in \mathbb{R}^{p \times T}$, corresponding to time segments of MEG signals, we jointly diagonalize the set of covariance matrices $C^i = \frac{1}{T}X_iX_i^\top$ (experiment \textbf{(c)}).
This practical task is in the spirit of~\cite{pham2001blind}.
Results are displayed in Figure~\ref{fig:results}

\subsection{Discussion}

In experiment (a), where the dataset is exactly jointly diagonalizable,
we observe the expected quadratic rate of convergence
for both the proposed algorithm and Pham's method.
We also observe that breaking the model (experiments (b) and (c)) makes convergence fall back to a linear rate.
As expected, the cost function does not go to $0$ in those cases.
We observe that the proposed quasi-Newton algorithm outperforms Pham's method by about an order of magnitude on each experiment.

\bibliographystyle{unsrt}
\bibliography{biblio}

\end{document}